\input lecproc.cmm

\def\Tn{{\bbbt}^n}

\def\T1{{\bbbs}^1}

\def\B{ B}
\def\bbbb{B}
\def\DD{\Delta_q}
\def\Rn{{\bbbr}^n}
\def\q{ q}
\def\qq{\overline q}
\def\p{p}

\def\z{\zeta^t}                                                 

\def\norm#1{\left\Vert #1\right\Vert}   
\contribution{ Periodic Orbits for Hamiltonian systems in Cotangent Bundles}
\contributionrunning{ Periodic Orbits}

\fonote{Partially supported by an NSF postdoctoral grant.}
\author {Christophe Gol\'e}
\address {IMS, SUNY at Stony Brook}
\abstract   
{ We prove the existence of at least $cl(M)$ periodic orbits
for certain time dependant Hamiltonian systems on the cotangent bundle of an
arbitrary compact manifold $M$. These Hamiltonians are not necessarily
convex but they satisfy a certain boundary condition given by a Riemannian
metric on $M$.
We discretize the variational problem
by decomposing the time 1 map into a product of ``symplectic twist maps''.
A second theorem deals with homotopically non trivial orbits in
manifolds of negative curvature.}

\titlea{0}{Introduction}
 	The celebrated theorem of Poincar\'e-Birkhoff states the existence of 
at least two fixed points for an area preserving map of the annulus $\T1 \times
[0,1]$ which ``twists'' the boundaries in opposite directions.

	In the 60's, Arnold proposed a generalization of this theorem for
a time 1 map $F$ of a time dependent Hamiltonian of $\Tn \times \B ^n$ (where 
$\B ^n$ is the closed ball in $\Rn$). While the Hamiltonian condition naturally
generalizes the preservation of area, a linking of the boundary of each
 fiber ( a sphere in $\Rn$) with its image by $F$ in the boundary of
 $\Tn \times \B ^n$ was
to replace the twist condition. Arnold [Ar1] conjectured that such a map
has at least as many fixed points as a real valued function $\Tn$ has critical
points.
 The philosophy was that fixed points for symplectic maps
should arise from Morse theory and not, say, from Lefshetz theory.

	Later, in [Ar2], he explained how fixed
points theorems on the annulus could be derived from theorems on
the 2-torus, by glueing carefully two annuli together 
(see also [Ch1]). He thus transformed
the problem to one of fixed points of symplectic maps on a {\bf compact}
symplectic manifold. This last conjecture, which 
asserts that the number of fixed point for the map is at least
equal to the minimum number of
 critical points of a real valued function on the manifold, is what
got to be well known as the Arnold conjecture.

 However, it is unclear whether 
the glueing construction can be done (symplectically) in higher dimensions. 
Even if it could, one would (if one could) have to use existing proofs of
the Arnold conjecture (e.g. [F2]), which we think are substantially harder 
than the
techniques we use here (and do not deal with homotopically non trivial orbits
as our Theorem 2 does).

	In 1982, Conley and Zehnder [CZ 1] gave a first proof 
of the Arnold conjecture for the torus $\bbbt^{2n}$.
In the same article, they also gave a direct proof of Arnold's 
original conjecture
on $\bbbt^n\times\B^n$.

 However, they were not able to use the linking of spheres in its
full generality. Their result remains crucial since it was the first 
non perturbational one in this direction.
The boundary condition that they used is expressed on the Hamiltonian in the
following way.
Letting $(\q,\p)$ be the coordinates on $\Tn \times \B ^n$ which is
endowed with the canonical symplectic structure $d \p \land d\q$, they set:

$$
H(\q,\p,t)= \langle A\p,\p\rangle + \langle b, p \rangle
\ \hbox{\rm for}\  \norm \p \geq K , \leqno (0.1)
$$
where $A=A^t$ is a non degenerate $n\times n$ matrix and $b\in \bbbr^n$
. This condition implies 
the linking of spheres
at the boundary.

	We propose here a version of this theorem on the cotangent bundle
of an arbitrary compact manifold. We also find, in a second theorem,
orbits of all free homotopy classes (and large enough period).

The bulk of this work was done as I was on a Postdoctoral position at the 
Forschungsinstitut f\"ur Mathematik, E.T.H. Z\"urich. I would like to
express my deep gratitude to Prof. Moser and Prof. Zehnder for inviting
me there. I had some invaluable discussions with them
as well as with  my companions
 Fredy K\"unzle, Boris Hasselblat, Frank Josellis, to whom
I extend my thanks.
I am very much endebted to 
 Patrice LeCalvez, whose work is
the starting point of mine.

Special thanks to Maciej Wojtkowski, Claude Viterbo, Misha Bialy,
 Leonid Polterovitch, Phil Boyland and Dusa McDuff
 for their specific help on this work.

Finally, were it not for the narrow mindedness of the French immigration
office, this work would have been joint with Augustin Banyaga. I dedicate
this work to him.

\titlea{1}{Results and basic ideas}

Let $(M,g)$ be a compact Riemannian manifold. Define
$$
B ^* M =\{(\q,\p) \in T^*M\ \ | \ \ g(\q)(\p,\p) = {\norm \p} ^2 \leq C^2
< R^2\},
$$
where $R$ is the radius of injectivity of $(M,g)$.
Let $\pi$ denote the canonical projection $\pi: B ^*M \rightarrow M$.

\theorem {1}{
Let $F: B ^*M \rightarrow B ^*M$ be the time 1 map of a time dependent
Hamiltonian $H$ on $B ^*M$, where $H$ satisfies the boundary condition:

$$
H(\q,\p,t)= g(\q)(\p,\p) \ \hbox{\rm for}\ \norm \p =C .
$$
Then $F$ has $cl(M)$ distinct fixed points and $sb(M)$ if they are all non 
degenerate. Moreover, these fixed points can all be chosen to to correspond
to homotopically trivial closed orbits of the Hamiltonian flow.}

	We remind the reader that $cl(M)$ is the cup length of $M$, which is
known to be a lower bound for the number of critical points of any real valued
function on $M$. Non degenerate means that no Floquet multiplier is equal to
one. $sb(M)$ is a lower bound for the number of critical points for a Morse
function on $M$.

\medskip {\bf Remark 1.2}{ It is important to note that, in the case where $M$ has $\bbbr^n$
as covering space, Theorem 1 can be expressed for a lift of $F$. In this case,
the radius of injectivity may be $\infty$ (e.g. for a metric close to a  flat
metric
on the torus, or when $M$ has a metric of negative curvature), and the set
 $B^*M$ can be as big as one wants. Theorem 1 can then serve as a starting
point to study Hamiltonian systems with asymptotic boundary conditions.}

\theorem{2}{Let $F$ be as in Theorem 1.
If (M,g) is of negative curvature, then $F$ has at least two
periodic orbits of period $d$ in  any given free homotopy class,
provided $d$ is big enough. In particular, $F$ has infinitely many
periodic orbits in $B^*M$.}
Exactly how big $d$ should be in Theorem 2
depends only on the metric. For a more precise 
statement, see section 7. Note also that if $H$ is 1--periodic in time,
periodic orbits of $F$ correspond to periodic orbits of the Hamiltonian
flow of $H$. Such a Hamiltonian system will then have infinitely
periodic solutions in $B^*M$.
\fonote{
If $H$ is not 1--periodic, periodic orbits of $F$ will correspond to
orbits of the Hamiltonian flow that
come back to their starting point, but generally at an angle. One can find
infinitely many of these orbits from Theorem 1, by applying it to
time $t$ map, $t\in (0,C)$, rescaling the metric each time.}

Note the difference in the boundary conditions (1.1)
and that of Conley-Zehnder
(0.1)
: theirs allow basically
all pseudo Riemannian metrics that are completely integrable and constant.
Ours only deals with Riemannian metrics, but with no further condition.
Note also that the orbits they find are homotopically trivial. We refer the
reader to [G1,2], [J] for the study of the homotopically nontrivial
case for $M=T^*\bbbt^n$ (the former with a method akin to that of this paper,
the latter in the spirit of [CZ 1]).

The method used to prove Theorem 1 and 2 is quite different from that of 
 Conley and 
Zehnder: whereas they use cut-offs on Fourier expansions,
we decompose the time 1 map into ``symplectic twist maps'' to get a finite
dimensional variational problem.

Symplectic twist maps are a natural
generalization of  monotone twist maps of the cotangent bundle of 
the circle (i.e. the annulus). 

In short, a symplectic twist map is a diffeomorphism $F$ from some
neighborhood $U$ of the zero section of $T^*M$ onto itself with the  property
that $ F^*pdq -pdq = dS$ for some $S$ and that
$(q,p)\to(q,Q)$ is a change of coordinates, where $F(q,p)= (Q,P)$.

To give an example, we make the following trivial remark. The 
shear map of the annulus :
$$
(q,p) \to (q+p,p),
$$
 which is a key model in the twist map
theory, is nothing more than the time 1 of the Hamiltonian $H_0(q,p)=
{1\over 2}p^2$ and in fact, its first coordinate map: 
$$
\eqalign{T^*\bbbs^1 & \to \bbbs^1\cr
(q,p) &\to q+p\cr}
$$
 is just the exponential
map for the standard (flat) metric on $\bbbs^1 \times \bbbr =T^*\bbbs^1$.

This suggest that the key model for symplectic twist maps on
the cotangent bundle $T^*M$ of a general compact manifold  $M$ should be
the time one map of a metric. The twist condition is given
in that case by the fact that the exponential map is a diffeomorphism 
of a neighborhood of zero in each fiber $T^*_qM$ and a neighborhood
of $q$ in $M$. Of course, most of the time, such a map is not completely
integrable.

If $F$ is a symplectic twist map, we have a simple proof of the original
 conjecture of Arnold:

\theorem {3}{(Banyaga, Gol\'e [BG]
) Let $F$ be a symplectic twist map of $B^*M$.
 Suppose that
each sphere $\partial B^*_qM$ links with its image by $F$ in 
$\partial B^*M$. Then the fixed points of $F$ are given	
by the critical points of a real valued function on $M$.}

In Appendix B, we reproduce the proof of [BG], for the convenience
of the reader. 
As in all these questions about fixed points, the major task is to make
the argument global: symplectic twist maps should be seen as local objects 
(even though
they should not be seen as perturbations) and the problem is to piece them 
together to form global ones. Here is one fundamental principle involved
in this.

Suppose we have two ``exact symplectic maps'':
$$
F^*pdq-pdq=dS \hbox{ and } G^*pdq -pdq= dS'
$$
Then it is simple to see that:
$$
(F\circ G)^*pdq -pdq= d(S\circ G +S')
$$
which we express as : generating functions add under compositions of maps.
This simple fact is key to the 
method in this paper: the functional
we use is a sum of generating functions of a finite sequence
of twist maps that decompose the time 1 map we study.

	This additivity property is
 the common thread between the method exposed here
and that of ``broken geodesics''
 reintroduced in symplectic geometry by Chaperon [Ch2]. The essential difference
is in the choice of coordinates in which one  expresses the generating function
:$(p,Q)$ in the method of Chaperon, $(q,Q)$ in the twist map method. In this
sense the twist map method is closer to the original method of broken 
geodesics as discribed in [Mi]. It  even coincides with it in the
case of the geodesic flow. 

Whereas Moser [Mo2] noticed that the time 1 map of a two dimensional
convex Hamiltonian can be decomposed into a product of twist maps,
the idea of decomposition of 
a time 1 map of a general Hamiltonian
stems from the  work of LeCalvez [L] on twist maps of 
the annulus. We generalize his simple but extremely efficient construction
to any cotangent bundle (Lemma (3.4)).

There are various theorems on the  suspension of  certain classes of
symplectic twist maps 
by Hamiltonian
flows ([D], [Mo2], [P-B]).
 In this sense, one might decide to forget about symplectic twist maps
and concentrate on Hamiltonian systems instead.
	In this paper, we take the opposite point of view: we think that
symplectic twist maps are a very useful tool to study Hamiltonian systems
 on cotangent bundles (see also the work of LeCalvez [L] on the torus).

The rest of the paper is organized as follows:

In section 2, we review some facts about geodesic flows and exponential maps.
We prove a lemma which is crucial for  the construction of an isolating block
in section 4.

In section 3 we give a precise definition of
 symplectic twist maps and prove the 
Decomposition Lemma (3.4).

In section 4, we use this decomposition and the additive property of generating
functions to construct a finite dimensional variational problem, 
i.e. a functional
$W$ on a finite dimensional space. This method is basically Aubry's ( [Au],
[Ka])
, when seen
on maps of the annulus.

In section 5, we construct an isolating block for the functional $W$.
For this, the boundary condition in the theorem is crucial.

In section 6, we make use of a theorem of Floer [F1] on global continuation 
of normally hyperbolic invariant sets: we
exhibit such an invariant set for the time 1 map of $H_0$ whose cohomology
survives under a deformation to our $H$. We then use the Conley-Zehnder Morse
theory to finish the proof of Theorem 1.

In section 7, we show how to adapt the proof of theorem 1 to the case
of non trivial homotopy classes, and prove theorem 2.

In Appendix A, we outline the connection that there is between the index
of the Hessian
of $W$ and the Floquet multipliers along a closed orbit of $F$.
 This is used in sections 5 and 7 to prove 
normal hyperbolicity of the invariant set.

 In Appendix B, we reproduce the proof of Theorem 3, given in 
[BG].

\titlea{1}{A few facts about the geodesic flow}
We start with some notation. 
Let $(M,g)$ be a Riemannian manifold. Both the tangent
fiber $T_qM$ and the cotangent fiber $T^*_qM$ are endowed with bilinear
forms:

$$\eqalign{&(v,v')\rightarrow g(q)(v,v ') \hbox{ for } v , v' 
\in T_qM,
\hbox{  and  }\cr
&(p,p ')\rightarrow g^\#(q)(p,p ') \hbox { for } p, p ' \in 
T^*_qM.\cr}$$
 
We will denote by $$\norm v :=\sqrt{g(q)(v,v')},\ { and }\ 
\norm p := \sqrt{g^\#(q)(p,p)},$$ 
hoping that the context will make it clear whether we speak about a vector or 
a covector.

The relation between $g$ and $g^\#$ is better understood in local coordinates:
If $A(q)$ denotes the matrix for $g^\#$ then $A^{-1}(q)$ is the matrix for
$g$. The terms of these matrices are usually denoted $g^{ij}$ and $g_{ij}$
respectively. The matrix $A(q)$ also gives the standard ( although $g$--
dependent) isomorphism between $T^*_qM$ and $T_qM$, which is an isometry 
for the above metrics. We will use the same notation ``$A(q)$''
for this isomorphism, even though it is coordinate independent, whereas the
matrix is not.

We want to outline here  some connections between the geodesic flow 
for
the metric $g$, the exponential map and the Hamiltonian flow for the Hamiltonian:

$$
H_0(q,p)= {1\over 2}{\norm p}^2.
$$
Let $T^*M$ be given the usual symplectic structure $dp \wedge dq$, and canonical
projection $\pi$. Let $h_0^t$ denote the time $t$ map of the Hamiltonian flow
of $H_0$. Then:
$$
exp_q(tA(q)p)=\pi\left(h^t_0(q,p)\right),
$$

This is basically a rewording of the equivalence of Hamilton's and 
Lagrange's equations under the Legendre transformation.
Here $H_0$ and $L_0(q,\dot q)={1\over 2} {\norm {\dot q}}^2$
 are Legendre transforms
of one another
under the change of coordinate $\dot q =\partial H_0 \bigm/
\partial p= A(q)p$ ([Arnold], section 15 or [Abraham-M] theorem 3.7.1 and 
3.6.2). This change of coordinate we will refer to as the Legendre 
transformation as well.

What is usually called the geodesic flow is just the flow $h_0^t$ restricted
to the (invariant) energy level $\{(q,p) \in T^*M \quad|\quad H_0(q,p)=1=
\norm p\}$ (the unit sphere bundle).

Because the exponential map :
	$$
\eqalign{exp&: TM \to M\times M\cr
	    &(q,v) \to (q,Q):=(q,exp_q(v))\cr}\leqno (2.1)
$$
defines a diffeomorphism between a neighborhood of the $0$--section in $TM$
and some neighborhood of the diagonal in $M\times M$ ([Mi], Lemma 10.3)
, we also have, via the 
Legendre transformation:
$$
\eqalign{exp&: T^*M \to M\times M\cr
            &(q,p) \to (q,Q):=\left(q,exp_q(A(q)p\right)\cr}
$$
which gives a diffeomorphism between a neighborhood of the $0$--section 
in $T^*M$ and some neighborhood of the diagonal in $M\times M$.
Just how big these neighborhoods are is measured by the radius of
injectivity $R$.

Because the Legendre transformation $A(q)$ is an isometry,
equation (2.1) gives a relation between  distances   between points in
$M$ and norms of vectors in $T^*M$:
$$
(q,Q)=exp(q,p) \Rightarrow \norm p = Dis(q,Q)
$$
It will be of interest for us to know the differential of the map
``Dis''.

\lemma{2.2}{If $(q,Q)=exp(q,p)$, and $h_0^1(q,p) = (Q,P)$, then:
$$
\partial_1Dis(q,Q)=- {p\over \norm p} \ \hbox{ and }\ 
\partial_2Dis(q,Q)= {P \over\norm P}
\leqno (2.2)
$$}

\proof { Let $v=A(q)p$. We have $Dis(q,Q)=\norm v$. The point 
$exp_q(-tv)$ is on the same geodesic as the one running from
$q$ to $Q$, namely $\{exp_q(tv)\bigm | t\in [0,1]\}$.
 For all small and positive $t$
we must then have:
$$
Dis(exp_q(-tv),Q)= (1+t)\norm v.
$$
Differentiating with respect to $t$ at $t=0$ yields:

$$
-\partial_1Dis(q,Q).v = \norm v .\leqno (2.3)
$$
On the other hand, by Gauss' lemma ([Mi], Lemma 10.5), the geodesic
through $q$ and $Q$ must be orthogonal to the sphere centered at $Q$ and
of radius $Dis(q,Q)$. This sphere is just the level set of the function:
$$
q' \to Dis(q',Q)
$$
whose gradient $A(q)\partial_1Dis(q,Q)$ at $q$ must be colinear to $v$. 
Equation
(2.3) yields:
$$
A(q)\partial_1Dis(q,Q)={-v\over \norm v}
$$
which immediately translates to the first equation we wanted to prove.

For the proof of the second equation, one must remember
that $V=A(Q)P$ is tangent at $Q$ 
to the geodesic between $q$ and $Q$ and has same
norm as $v$. It is , more precisely, the parallel transport of $v$ along
this geodesic. Thus:
$$
\bigl.{d\over dt}Dis\left(q,exp_q\left((1+t)v\right)\right)\bigr|_{t=0}
=
\norm V = \partial_2 Dis(q,Q).V
$$
and the rest of the reasoning is the same as for the first equation.
\qed}

\titlea{3}{Symplectic twist maps and the decomposition lemma}

If $H(q,p,t)$ is an optical Hamiltonian function (i.e. $H_{pp}$ is convex), then
its flow has many similar features to that of 
$H_0(q,p)={1\over2}{\norm p}^2$. In particular if $F$ is its time $\epsilon$
, and
$F(q,p)=(Q,P)$, the correspondance $(q,p) \to (q,Q)$ is a diffeomorphism 
between suitable neighborhoods of the 0--section in $T^*M$ and the diagonal
in $M \times M$ (compare equation 2.1).
 This can be seen in a chart, looking at the Taylor series of the 
solution with respect to small time:

$$
\eqalign{Q=q(\epsilon)&= q(0)+ \epsilon.H_p+ o(\epsilon^2)\cr
	 P=p(\epsilon)&= p(0) -\epsilon.H_q + o(\epsilon^2)\cr},
$$

>From this we see that ${\partial Q \over\partial p}(z(0))$ 
is non degenerate. This remark was made by Moser in the dimension 2
case ([Mo2]).

Another feature enjoyed by Hamiltonian flows is that they are
 exact symplectic.

 These two properties put together give us
the following:

\definition {3.1}{A {\bf symplectic twist map} F is a diffeomorphism
of a neighborhood $U$ of the 0--section of $T^*M$ onto itself satisfying
the following:
\medskip
\item{(1)}$ F$ is {\bf exact symplectic}:\ $ F^*pdq-pdq=dS$ for some real 
function $S$
on 
$U$.
\item{(2)} ({\bf Twist})
if $F(q,p)=(Q,P)$, then the map $\psi:(q,p) \to (q,Q)$ is
	embedding of U in $M\times M$.
\medskip 
The function $S(q,Q)$ is then called the {\bf generating function} for $F$.
}
\medskip
 {\bf Remark 3.2}{Of course ([G1,2,3]), monotone twist maps of the annulus 
(i.e. of $T^*\bbbs ^1)$) are
symplectic twist maps in the sense of this definition. $U$ is usually taken
to be either the whole cylinder, or the subset $\bbbs^1 \times [0,1]$.
Note that one way to express the twist condition is by saying that the image
by $F$ of a (vertical) fiber in $U$ intersects any fiber in at most one
point.

 To my knowledge, the term symplectic twist map was introduced by
McKay, Meiss and Stark. Their definition ([MMS]) is a little more restrictive
than the above, in that they work on $T^*\bbbt^n$ and ask that $\partial Q/
\partial p$ be definite positive. Our condition only implies that $det (
\partial Q/\partial p)\not= 0$. Similar maps have also been studied
extensively by Herman ([H]): he called them monotone. We also have used this 
terminology ([G1,2]) but in the end found it misleading because we were also
dealing with monotone flows [G2], the two concepts being only 
related in certain cases.}
\medskip 
{\bf Remark 3.3}
{Equation (2.1) tells us that the time 1 map $h^1_0$
of $H_0$ is also a symplectic twist map on some neighborhood $U$ of 
the 0--section. 
 Note that for the time 1 of an
 Hamiltonian H, the function S is (when defined) the 
action:
$$
S(q,Q)=\int_{(p,q)}^{(P,Q)}pdq-Hdt
$$
taken along the unique solution of the Hamiltonian
flow between $(p,q)$ and $(P,Q)$.
If $L$ is the Legendre transform of $H$, the above integral is just 
$$
S(q,Q)=\int_0^1L(q,\dot q,t)dt
$$
along the solution. In  the case where $H=H_0$, $L(q,\dot q,t)={1\over 2}
\norm {\dot q}^2$, i.e. $S$ is the energy of the (unique) geodesic between
$q$ and $Q$.

As noted in the introduction, $h_0^1$ should be
our model map, the way the shear map is the model map in the theory of
monotone twist maps.}

The reason why twist maps can be so useful lies
in  the following fundamental lemma, due to LeCalvez [L]
in the case of diffeomorphisms of the annulus isotopic to the $Id$:

\lemma{3.4 (Decomposition)}{ (LeCalvez, Banyaga, Gol\'e): 

 Let $F$ be the time 1 map of a (time dependent) Hamiltonian on
a compact neighborhood $U$ of the 0--section. Suppose that $F$ leaves
$U$ invariant. Then, $F$ can be
decomposed into a finite product of symplectic twist maps:
$$
F= F_{2N} \circ \ldots \circ F_1
$$
}
\medskip 
{\bf Remark 3.5}No convexity is assumed of the Hamiltonian, nor any closeness
to an integrable one.

\proof { Let $G_t$ be the time $t$ map of the Hamiltonian, 
starting at $t=0$. We can write :
$$
F= G_1\circ G_{N-1 \over N}^{-1}\circ \ldots \circ G_{k\over N}\circ 
G_{k-1\over N}^{-1}\circ \ldots \circ G_{1\over N}\circ Id
$$
and each of the $G_{k\over N}\circ G^{-1}_{k-1\over N}$ is an exact symplectic map,
which we can make as close as we want to the $Id$ by increasing $N$.
If $H_{pp}$ is positive definite, each of these maps are twist, by Moser's
remark, and we are done ($F$ is the product of $N$ twist maps in this case)
. In general, we do the following.
The twist condition (2) in Definition (3.1) of symplectic twist is an open 
condition.
Hence, if $h_0^t$ is the time $t$ map of $H_0$,
 the map $F_{2k-1}:= h_0^{-1}\circ G_{k\over N}\circ G^{-1}_{k-1\over N}$ must
satisfy (2) for $N$ big enough (here, the compactness of $U$ is needed).
We then set $F_{2k}=h_0^1$ for all $k$ to get the decomposition advertized.\qed}

\medskip
{\bf Remark 3.6}{We leave it to the reader to check that
Lemma 3.4 is also valid for lifts of maps to the covering space of $M$.}

\titlea{4}{The discrete variational setting}

Let $F$ be as in Theorem 1. From the previous section, we can write
$$
F= F_{2N} \circ \ldots \circ F_1,
$$
with the further information that $F_{2k}$ restrained to the boundary 
$\partial B^*M$ of $B^*M$ is the time 1 map of $H_0$, that we have called
$h_0^1$. Likewise, $F_{2k-1}$ is $h_0^{{1\over N}-1}$ on $\partial B^*M$,
by the proof of 
the decomposition Lemma (3.4) and the boundary condition (1.1) imposed on $F$.

	Let $S_k$ be the generating function for the twist map $F_k$ and
$\psi_k$ the diffeomorphism $(q,p) \to (q,Q)$ induced by $F_k$. We can 
assume that $\psi_k$ is defined on a neighborhood $U$ of $B^*M$ in $T^*M$.
	Let 
$$
\eqalign{
O=\{\overline q =(q_0,\ldots,q_{2N-1})\in  M^{2N}\mid &(q_k,q_{k+1})
\in \psi_k (U)  \hbox{ and } \cr
&\quad (q_{2N},q_0) \in \psi_{2N-1}(U) \}\cr}
\leqno (4.1)
$$
 $O$ is an open set in $M^{2N}$, containing a copy of
$M$ (the elements $\overline q$ such that $q_k=q_0$, for all $k$).

Next, define :
$$
W(\overline q)= \sum _{q=0}^{2N-1} S_k(q_k,q_{k+1}), \leqno (4.2)
$$
where we have set $q_{2N}=q_0$.
Let $p_k$ be such that $\psi_k(q_k,p_k)=(q_k,q_{k+1})$ and
let $P_k$ be such that $F_k(q_k,p_k)= (q_{k+1},P_k)$. $P_k$ and $p_k$
are well defined functions of $(q_k, q_{k+1})$.

We claim:
\lemma{4.3}{The sequence $\overline q$ of $O$ is a critical point of $W$
if and only if the sequence $\{(q_k,p_k)\}_{k\in \{0,\ldots,2N,0\}}$ is an
orbit under the successive $F_k$'s, that is if and only if $(q_0,p_0)$
is a fixed point for $F$.}

\proof{Because the twist maps are exact symplectic and using 
the definitions
of $p_k$, $P_k$, we have:
$$
P_kdq_{k+1}-p_kdq_k=dS_k(q_k,q_{k+1}),\leqno (4.4)
$$
and hence
$$
dW(\overline q)= \sum_{k=0}^{2N-1}(P_{k-1}-p_k)dq_k
$$
which is null exactly when $P_{k-1}=p_k$, i.e. when 
$F_k(q_{k-1}, p_{k-1})= (q_k, p_k)$. Now remember that we assumed
that $q_{2N}=q_0$.
\qed }

Hence, to prove Theorem 1, we need to find enough critical points for $W$.
For this , we will study the gradient flow of $W$ (where the gradient will
be given in terms of the metric $g$) and use the boundary condition to find
an isolating block. 

We now indicate how this variational setting is related to the classical
method of broken geodesics, and how to modify it to deal with homotopically
non trivial solutions.

Because each $F_k$ is close to $h_0^{t_k}$ for some positive or negative
$t_k$,  we have that:
$$
q \in \psi_k (B^*M_{q})
$$
and, since $B^*_qM \to \psi_k(B^*M)$ is a diffeomorphism, we can define
a path
$c_k(q,Q)$ between $q$ and a point $Q$ of $\psi_k(B^*_qM)$
by taking the image of the oriented line segment between $\psi_k^{-1}(q)$ and
$\psi_k ^{-1}(Q)$ in $B^*_qM$. In the case where $F_k=h_0^1$, this amounts
to taking the unique geodesic between $q$ and $Q$ in $\psi_k(B^*_qM)$
.            %

If we look for periodic orbits of period $d$ and of a given homotopy type,
we decompose $F^d$ into $2Nd$ twist maps, by decomposing $F$ into $2N$.
Analogously to (4.1), we define :

$$
\eqalign{
O_d=\{\overline q =(q_0,\ldots,q_{2Nd-1})\in  M^{2Nd}\mid &(q_k,q_{k+1})
\in \psi_k (U)  \hbox{ and }\cr
& (q_{2Nd},q_0) \in \psi_{2Nd-1}(U)\} ,\cr}
$$
remarking that the $\psi_k$'s here correspond to the decomposition of $F^d$
 into $2Nd$ steps ($U$ is as before a neighborhood of $B^*M$).

To each element $\qq$ in $O_d$, we can associate a closed curve, made by
joining up each pair $(q_k,q_{k+1})$ by the unique curve
$c_k(q_k,q_{k+1})$ defined above. This loop $c(\qq)$ is piecewise 
differentiable
and it depends continuously on $\qq$, and so does its derivatives (left
and right). In the case of the decomposition of $h_0^1$ , taking
$F_k$=$h_0^1$, this is exactly the construction of the broken 
geodesics ([Mi], \S 16). Now any closed curve in $M$ belongs to a 
free homotopy class $m$.

To any $d$ periodic point for $F$, we can associate a sequence $\qq (x)
 \in O_d$
of $q$ coordinates of the orbit of this point under the successive $F_k$'s
in the decomposition of $F^d$.

\definition{4.5}{Let $x$ be a periodic point of period $d$ for $F$. Let $\qq$ 
be the sequence in $O_d$ corresponding to $x$. We say that $x$ is an $(m,d)$ 
point if $c(\qq (x))$ is in the free homotopy class $m$.}

To look for $(m,d)$ orbits in (Theorem 2 in section 7), we will work in:

$$
O_{m,d}=\{\qq \in O \mid c(\qq) \in m \} \leqno (4.6)
$$
Since $c(\qq)$ depends continuously on $\qq \in O$, we see that 
$O_{m,d}$ is actually a connected component of $O$.

The functional $W$ will be given this time by:
$$
W(\qq)=\sum_{k=0}^{2Nd-1}S_k(q_k,q_{k+1})
$$ 
defined on $O_{m,d}$.
Again, as in Lemma 4.3, critical points of $W$ in $O_{m,d}$ correspond
to $(m,d)$ periodic points.

\medskip 
{\bf Remark 4.7}{
The reader that wants to make sure that, in the proof of Theorem 1,
the orbits found are homotopically trivial, should check that
throughout the proof, one can work in the component
$O_{e,1}$ of $O_1=O$ of sequences $\qq$ 
which have $c(\qq)\in e$, where $e$ is the identity element of $\pi_1(M)$.
}

\titlea{5}{The isolating block}

In this section we prove that the set $B$ defined as follows:
$$
B=\{\overline q \in O\mid \Vert p_k(q_k,q_{k+1}) \Vert \leq C\} \leqno (5.1)
$$
is an isolating block for the gradient flow of
$W$, where $O$ is defined in (4.1),
 $C$ is as in (1.1) and $p_k$ is the function defined
in the previous section (see below (4.2))
. To try to visualize this set in $M^{2N}$, the reader 
should realize that the twist condition on $F_k$ and the fact that $F_k$
coincides with the time 1 or time ${1\over n}-1$ of the hamiltonian $H_0$ at
the boundary of $B^*M$ implies that:
$$
 Dis (q_k,q_{k+1})= a_k\norm {p_k} 
\quad\hbox{ where } \left\{\matrix{a_k&=&1& \hbox{ if $k$  is even} \cr
				  a_k&=&{1-N\over N}&\hbox{if $k$ is odd} \cr}
\right.\leqno (5.2)
$$

Note that $B$ still contains a copy of $M$ ( the constant sequences).

We will define an {\bf isolating block} for a flow  to be a compact
neighborhood with the property that the solution through each boundary point
of the block goes immediately out of the block in one or the other time 
direction ( [C], 3.2 ). Sometimes, more refined definitions are made,
but this one is sufficient to ensure that the maximal invariant set for
the flow contained in the block is actually contained in its interior:
a block in this sense is an isolating neighborhood, which is really the only 
property we need here.

\proposition{5.3}{ $B$ is an isolating block for the gradient flow of $W$.}

\proof{ Suppose the point $\overline q$ of $U$ is in the boundary of $B$.
this means that $\norm p_k =C$ for at least one $k$. As noted in (5.2), this
means that $Dis (q_k, q_{k+1})= a_kC$ for some factor $a_k$ only depending on
the parity of $k$. We want to show that this distance  increases either
in positive or negative time under the gradient flow of $W$. This flow
is given by:
$$
\dot q_k = A_k(P_{k-1}-p_k)=\nabla W_k(\qq)
\leqno (5.4)
$$
Where $A_k=A(q_k)$ is the duality morphism associated to the metric $g$ at
the point $q_k$ (see beginning of section 2). Remember that we have put
the product metric on $O$, induced by its inclusion in $M^{2N}$.
 
Let us compute the derivative of the distance along the flow at a boundary 
point of $B$, using Lemma 2.1:
$$
\eqalign{\bigl.{d\over dt} Dis(q_k,q_{k+1})\bigr|_{t=0}
&= \partial _1 Dis (q_k,q_{k+1})
.\nabla W_k(\overline q)\cr
&\qquad +\partial _2 Dis (q_k,q_{k+1}). \nabla W_{k+1}
(\overline q)\cr
&= ({a_k\over |a_k|}){-p_k\over \norm {p_k}}.A_k.(P_{k-1}-p_k) 
\cr 
& \qquad+({a_{k}\over|a_{k}|}){P_k\over \norm {P_k}}.A_{k+1}.(P_{k}-p_{k+1})
\cr}
\leqno (5.5)
$$
We now need a simple linear algebra lemma to treat this equation.

\lemma{5.6}{Let $\langle \, , \,\rangle$ denote a metric form in $\bbbr^n$, 
and
$\norm {. }$ its corresponding norm. Suppose that $p$ and $p'$ are in $\bbbr^n$
,that 
$\norm p=C$ and  that
$\norm { p'} \leq C$. Then :
$$
\langle p\ ,\ p'-p\ \rangle \leq 0.
$$
Moreover, equality occurs if and only if $p'=p$.
}
\proof {From the positive definiteness of the metric, we get:
$$
\langle \ p'-p, p'-p\ \rangle \geq 0,
$$
with equality occuring if and only if $p'=p$ (call this last
assertion *). From this, we get:
$$
2\langle \ p, p'\ \rangle \leq \langle \ p', p'\ \rangle +
\langle\ p, p\ \rangle
$$
with *.
Finally,
$$
\langle \ (p'-p),p \ \rangle = \langle \ p',p \ \rangle - \langle \ p,p
\  \rangle \leq 0
$$
with *.
\qed}
Applying  Lemma 5.6 to each of the right hand side terms in (5.5),
we can deduce that ${d\over dt}Dis(q_k, q_{k+1})$ is
positive when $k$ is pair, negative when $k$ is odd. Indeed,
because of the boundary condition in the hypothesis of the theorem,
we have $\norm {P_{k}}= \norm {p_{k}}$ whenever $\norm {p_k} =C$:
the boundary $\partial B^*M$ is invariant under $F$ and all the $F_k$'s.
On the other hand $\overline q \in B \Rightarrow \norm {p_l}\leq C
\hbox{ and }
 \norm {P_l}
 \leq C$, for all $C$, by invariance of $B^*M$
. Finally, $a_k$ is positive when $k$ is even, negative when
$k$ is odd.

But what we really want is this derivative to be of a definite sign, not
zero. It is certainly the case when at least one of $\nabla W_k(\qq ), 
\nabla W_{k+1}(\qq)$
is not zero. Suppose they are both zero. Then $k$ is in an  interval
$\{ l,\ldots ,m\}$ such that, for all $j$ in this interval, $\norm {p_j} =C
=\norm {P_j}$ and $\nabla W_j(\qq) =0$.

 It is now crucial to notice that 
$\{l, \ldots ,m \}$ can not cover all of $\{0,\dots,2N\}$: this would mean
that $\qq$ is a critical point corresponding to a fixed point of $h_0^1$
in $\partial B^*M$. But such a fixed point is forbidden by our choice of $C$:
geodesics in that energy level can not be fixed loops $(C>0)$, and they
can not close up in time one either ($C$ is less than the injectivity radius).

We now let $k=m$ in (5.5) and see that the flow must definitely
escape the set $P$ at $\qq$ in either
positive or negative time, from the the $m^{th}$ face of $P$.
\qed}
\medskip
{\bf Remark 5.7}{If we have decomposed the time 1 map of a Hamiltonian that is 
positive definite into a 
product of  $N$ twist maps, all the $F_k$'s coincide 
with $h_0^{1\over N}$ on the boundary of $B^*M$. In that case,
$$
\norm {p_k} = {1\over N}Dis (q_k,q_{k+1}), \hbox{ for all } k
$$
and the $a_k$'s in the above proof are always positive. Following the argument
through, we find that $B$ is a repeller block in this case: all points on 
$\partial B$ exit in positive time.}
\medskip
{\bf Remark 5.8}{LeCalvez ([L]) provides a more detailed analysis of the behavior
of the flow at ``corner'' points of his analog of the set $B$. He indicates
an induction to show that the flow enters or exits the $j^{th}$ face 
($j$ is in $\{l,\ldots,m\}$ as in the above proof) at different orders in small
time. Such a reasoning could be made in our context also, but we find it
unnecessary, given our working definition of an isolating block.}

\titlea{6}{Proof of Theorem 1}

To finish the proof of Theorem 1 we will be using a refinement of the
Conley Index continuation proved by Floer ([F1]). The homology
group of the invariant set $G^\lambda$ appearing in
this lemma
bears the germs of what became
later  Floer Cohomology (see e.g. [F2], and also [McD])
, and in the case that we study, it is
probable that it is one and the same thing. The present approach enables us
to avoid the problem of infinite dimensionality
in [F2], i.e. all the analysis.

\lemma{6.1} {(Floer) Let $\phi^t_\lambda$
 be a one parameter family of flows on a $ C^2 $ manifold
$M$. Suppose that $G^0$ is a compact $C^2$ submanifold invariant under the flow
$\phi^t_0$. Assume moreover that $G^0$ is {\bf normally hyperbolic}, i.e. there is a 
decomposition:
$$
TM_{|G^0}= TG^0 \oplus E^+ \oplus E^-
$$
which is invariant under the covariant linearization of the vector field
$V_0$ corresponding to $\phi^t_0$
with respect to some metric $\langle\ ,\ \rangle$, so that for some
constant $m>0$:
$$
\eqalign{
\langle\xi, DV_0 \xi\rangle &\leq -m\langle\xi,\xi\rangle \hbox{ for }
\xi \in E^-\cr
\langle\xi, DV_0 \xi\rangle &\geq m\langle\xi,\xi\rangle \hbox{ for }
\xi \in E^+\cr}
\leqno{6.1}
$$
Suppose that there is a retraction $\alpha : M\to G^0$ and that there
is a compact neighborhood $B$ which is isolating for all $\lambda$.
Then, if $G^\lambda$ denotes the maximum invariant set for
$\phi_\lambda^t$ in $B$, the map:

$$
\bigl(\alpha_{|G^{\lambda}}\bigr)^*: H^*(G^0) \to H^*(G^\lambda)
$$
in \v Cech cohomology
is injective.}

In this precise sense, normally hyperbolic invariant sets continue
{\it globally}: their topology can only get more complicated as the
parameter varies away from $0$. Note that we have given here a
watered down version of Floer's theorem. His uses the notion of
Conley continuation of invariant sets. He also works in the equivariant
case. But the above, taken from  his Theorem 2 in [F1], is what we need here.

 The family of flows we consider is $\z _\lambda$, the flow solution
of 
$${d\over dt}\overline q = \nabla W^\lambda (\qq),$$
 and $W^\lambda$
is defined as in 4.2 for the map $F_\lambda$, time 1 map of the
Hamiltonian:
$$
H_\lambda = (1-\lambda)H_0 + \lambda H
$$
We can assume that this construction is well defined, i.e., that we
make the decomposition in the Decomposition Lemma 3.4 fine enough to fit
any $F_\lambda, \lambda$ in $[0,1]$.
The manifold on which we consider these (local) flows is $O$, an open 
neighborhood
of $B$ in $M^{2N}$.
 Of course, each of the $F_\lambda$
satisfies the hypothesis of Theorem 1, and thus Proposition 5.3
applies to $\z_\lambda$ for all $\lambda$ in $[0,1]$: $B$ is an isolating block
for each one of these flows.

	The part of Floer's lemma that we are missing so far is the normally
hyperbolic invariant manifold for $\z_0$.

\lemma{6.2}{Let $G^0= \{ \overline q \in B \mid q_k=q_0,  \forall k \}$.
Then $G^0$ is a normally hyperbolic invariant set for $\z_0$. It is a retract
 of $O$ and is the maximal invariant set in $B$.}

\proof{
All the $F_k$'s in the decomposition of $h^1_0$ are time $a_k$ maps of
the Hamiltonian $H_0$, for  $a_k$ as in (5.2). But for this Hamiltonian, the
$0$--section of $T^*M$ is made out of fixed points. These translate, in
terms of sequences, to points in $G^0$. Moreover, these are the {\it only}
periodic orbits for the Hamiltonian flow of $H_0$ in $B^*M$, by the definition
of this set. (e.g. in the case $M=\bbbs ^n$ with 
the standard metric, the orbits corresponding to great circles would not be 
fixed points of $h_0^1$ in $B^*M$).
	
This implies that $G^0$ is the maximum invariant set for $\z_0$ in $B$.
Indeed, since $\z_0$ is a gradient flow,
 such an invariant set should be formed by critical points and
connections between them. We saw that there are no other critical points 
but the points of $G^0$. If there were a connection orbit entirely lying
in $B$, it would have to connect two points in $G^0$, which is absurd
since by continuity any two points of $G^0$ give the same value for $W^0$,
whereas $W^0$ should increase along non constant orbits.

$G^0$ is a retract of $M^{2N}$ under the composition of the maps:
$$
\overline q=(q_1,\ldots ,q_{2N})\to q_1 \to (q_1,q_1,\ldots ,q_1)=\alpha
(\overline q)
$$
which is obviously continuous and fixes the points of $G^0$.

We are left to show that $G^0$ is normally hyperbolic. For this, we are going
to appeal to a relationship between the linearized flow of $\z_\lambda$ and
that
of $H_\lambda$. The following lemma was proven by McKay and Meiss in the
twist map of the annulus case. We present their proof in Appendix A: it
holds in the setting of general cotangent bundles.

\lemma{6.3}{( [M-M]) Let $F$ be the time 1 map of a Hamiltonian and let
$W$ be its associated functional.
If $\qq$ is a critical point corresponding to the orbit
of $(q_0,p_0)$, the set of eigenvectors of eigenvalue 1 of 
 $DF_{(q_0,p_0)}$ are in 1--1 correspondence with the set of
eigenvectors of eigenvalue 0 of $Hess W(\qq)$ }

To use this lemma,
we  remark that since $G^0$ is made out of critical points, saying
that it is normally hyperbolic is equivalent to saying that
$Hess W^0(\qq)$ has exactly $n=dim G^0$ eigenvalues equal to zero for any
point $\qq$ in $G^0$. These eigenvalues have to correspond to eigenvectors
in $TG^0$, the normal space of which must be spanned by eigenvectors with
non zero eigenvalues ($Hess W^0$ is symmetric). Hence, from 
 Lemma 6.3, it is enough to check that at a point $(q_0,0) \in B^*M$
corresponding to $\qq$, $1$ is an eigenvalue of multiplicity
exactly $n$ for  $Dh_0^1(q_0,0)$.
Let us compute $Dh_0^1(q_0,0)$ in local coordinates. It is the solution at time 1
of
the linearized equation:
$$
\dot U=JHess H_0 (q_0,0)U
$$
 along the constant solution $(q(t),p(t))=(q_0,0)$, where $J$ denotes the usual
symplectic matrix $\pmatrix{0&I\cr -I&0\cr}$.	
An operator solution for the above equation is given by
$
exp \left(t JHess H_0 (q_0,0)\right)
$
On the other hand:
$$
Hess H_0 (q_0,0)= \pmatrix{
0&0\cr
0&A(q_0)\cr}
$$
which we computed from $H_0(q,p)=A(q)p.p$, the zero terms appearing 
at $p=0$ because
they are either quadratic or linear in $p$.
>from this,
$$
Dh_0^1(q_0,0)= exp \left(J Hess H_0(q_0,0)\right)= \pmatrix{ I & A(q_0)\cr
						 0 & I }
$$
is easily derived. This matrix has exactly $n$ eigenvectors of eigenvalue 1
( it has in fact no other eigenvector). Hence, from Lemma 6.3, $Hess W(\qq)$
has exactly $n$ vectors with eigenvalue $0$, as was to be shown.\qed
}
We now conclude the proof of Theorem 6.3. We have proved that the flow $\z$,
which is gradient, has an invariant set $G=G^1$ with $H^*(M) \hookrightarrow
H*(G)$. From this we get in particular: 
$$
cl(G) \geq cl(M) \hbox{ and } sb(G) \geq sb(M).
$$
The corollary of Theorem 5 in [CZ 1] tells us that $\z$ must have at least
$cl(G)$ rest points in the set $G$, whereas the generalized Morse inequalities
in Theorem 3.3 of [CZ 2] tell us that, if they are all assumed to be non 
degenerate, $\z$ must have $sb(G)$ rest points. But Lemma 6.3 tells us 
that non degeneracy for $Hess W$ at a critical point is the same thing as 
nondegeneracy of a fixed point for $F$ (no eigenvector of eigenvalue 1).

As was stated in Remark 4.7, we could have worked in $O_{e,1}$ all along
to guarantee that the orbits found are homotopically trivial.
The only thing that one should check is that $G^0_{m,d}$ is indeed in this
component of $O$, which is the case.
 This concludes the proof of Theorem 1. \qed

\titlea{7}{Negative curvature and orbits of different homotopy types}

\titleb{7.1}{Setting the problem}
We are going to sketch here the changes needed in the proof of Theorem 1
in order to prove Theorem 2 on the existence of orbits of different
homotopy types.

It is a known (see e.g. [GHL], 2.98) that on a compact Riemannian
manifold there exists in any nontrivial free homotopy class $m$ a smooth
and closed geodesic which is of length minimal in $m$.
\fonote{
We remind the reader that free homotopy classes of loops differ from
elements of $\pi_1(M)$ in that no base point is kept fixed under the
homotopies. As a result, free homotopy classes can be seen as conjugacy
classes in $\pi _1(M)$, and thus can not be endowed with a natural algebraic
structure. Two elements of a free class give the same element in $H_1(M)$.
Hence free homotopy classes form a set smaller than $\pi_1(M)$, bigger
than $H_1(M)$. All these sets  coincide if $\pi_1(M)$ is abelian.}

	Moreover, a theorem of Cartan asserts that
if the manifold is of negative curvature, there is in fact
one and only one geodesic in each class $m$ ($m$ 
not containing the point curves
 ([Kl], Theorem 3.8.14)).

Let $M$ be of negative curvature and 	
let $l(m)$ denote the length of the geodesic in $m$ in that case.

In 4.5, we have defined $(m,d)$ orbits by saying that a certain curve
that the orbit defines in $M$ is of class $m$. We could 
also use a favored lift $\overline F$ of $F$ 
to the covering space $\overline M$ of $M$ to define
such orbits, by asking:
$$
\overline F^d(x)=m.x
$$
where $m.x$ denotes the action of $m$ seen as a deck transformation in
$T^*\overline M$
(the favored lift is the one corresponding to lifting the solution curves
of the Hamiltonian flow). But definition 4.5 turns out to be more convenient 
to use here (both are equivalent, of course).
We now restate:

\theorem{2} {Let (M,g) be a Riemannian manifold of negative curvature.
and $H$ be as in Theorem 1. Then, whenever $(M,g)$ has a geodesic
whose class in $\pi_1(M)$ is $m$, $F$ has at least 2 $(m,d)$ orbits in $B^*M$
when $l(m) <dC$ .}

\medskip 
{\bf Remark 7.1.1}{The fact that we find two orbits and not a number given by
the topology of the manifold is not an artifact of the
proof, but derives from the unicity of the closed geodesic in a given
class. Note also that we do not guarantee that an orbit of the form
$(m^k,kd)$ is not actually an $(m,d)$
orbit. We should then ask for $(m,d)$ to be
prime, in that very sense.

Note also, there are a priori more of these pairs $(m,d)$ than there are rational
homology directions.}

 The proof of Theorem 2 has the same broad outline as that of Theorem 1.

We decompose $ F=F_{2N}\circ \ldots \circ F_0$ as before, which
gives us a decomposition of $F^d$ into $2Nd$ twist maps.

We would like to claim, in analogy to Proposition 5.3 that
$$
B=\{ \qq \in O_{m,d} \mid \norm{p_k(q_k,p_k)} \leq C\}
$$
is an isolating block for the gradient flow of $W$. 

But this will not be enough for our purpose.
To make sure that two critical points correspond 
to points that are actually on 2 distinct
orbits, one should do the following: to decompose $F^d$, we have
decomposed $F$ in $2N$ steps. Define:
$$
\sigma: O_{m,d}\to O_{m,d} \hbox{ by setting }\  (\sigma \qq)_k=q_{k+2N}
$$
where we identify: $q_{k+2Nd}=q_k$.
It is clear that 2 critical sequences corresponding to points in the same 
orbit by $F$ get identified in the quotient by the action of $\sigma$. 
 So our candidate for isolating block will be given by the quotient
$B/\sigma$ of $B$ by the action of $\sigma$ (note that $\sigma$ leaves
$B$ invariant, so that the quotient makes sense).

It can be seen that if $m$ is non trivial, then the action of $\sigma$
is without fixed points. Since it is also periodic, the action is then
properly discontinuous ([Gr], Chapter 5), and hence the quotient map
$O_{m,d}\to O_{m,d}/\sigma$ is a covering map.

We now describe the candidate for normally hyperbolic invariant set.
It will be the quotient by $\sigma$
of the set $G^0_{m,d}$ made of the  critical sequences corresponding to
the continuum of $(m,d)$ orbits that form the closed geodesic of class
$m$, parametrized so that the Hamiltonian flow goes through it in time $d$.
Call this orbit $\gamma$.

Note that $G^0_{m,d}$
contains all the possible critical points
for $W^0$ in $ O_{m,d}$ since a critical point for $W^0$
must be contained in a continuum of critical points: if $(q_0,p_0)$
is an $(m,d)$ point, so is $h^t_0(q_0,p_0)$, for any $t$. But we know,
in the case of $H_0$ that there is one and only one such set, namely
$\gamma$.

Writing $F^0_k\circ \ldots \circ F^0_1= \phi_k$, where the $F^0_k$'s decompose
the map $h_0^d$, we can write:
$$
G^0_{m,d}=\{\qq (t) \in O_{m,d} \mid q_k =\pi\circ\phi_k(\gamma(t))\}
$$
where $\pi(q,p)=q$ is the canonical projection. Again, since $\sigma$
restricted to $G^0_{m,d}$ actually corresponds to the action of $F$ on 
$\gamma$, $G^0_{m,d}$ is $\sigma$ invariant and hence the quotient
$G^0_{m,d}/\sigma$ makes sense. Since the quotient map is a covering map
and $G^0_{m,d}\cong \gamma$, we have :
$$
G^0_{m,d}/\sigma \cong \bbbs^1 \leqno (7.1.2)
$$

\titleb{7.2}{Proof of Theorem 2}

\lemma{7.2.1}{$B/\sigma$ is an isolating block}
\proof{ Because we have assumed 
$l(m)<dC$, 
 $F$ cannot have any $(m,d)$ orbits confined to $\partial B^*M$: 
since $F$ coincides with $F_0$ on this set,
such
an orbit would have to correspond to a closed geodesic of
 free homotopy class $m$, but
of length $dC$, which is absurd.
This in turn implies that $W$ has no critical points on $\partial B$,
and the reasoning of Proposition 5.3 applies without change to show that
points in $\partial B$ must exit $B$ in positive or negative time.
Since the covering map is a local diffeomorphism, this is also true
in $B/\sigma$, which is then an isolating block. \qed }

\lemma{7.2.2}{$G^0_{m,d}/\sigma$ is a normally hyperbolic invariant set for
$\z_0$. It is a retract of $O_{m,d}/\sigma$.}

\proof{ We prove the statement ``upstairs'', taking the quotient by
$\sigma$ only at the end.

According to Lemma 6.3, and the reasoning in 
the proof of Lemma 6.2, it is enough to show that the differential of
$h_0^d$ on a point of  $\gamma$  has no other eigenvector 
of eigenvalue 1 than the 
vector tangent to $\gamma$.

To compute the differential of $h_0^d$ at the point $(q_0,p_0)=\gamma(0)$ 
 we are going to choose a coordinate
system $(z,t,s)$ around $(q_0,p_0)$ in the following way:
$z,t$ will be a coordinate system for a tubular neighborhood in the
energy surface containing 
$(q_0,p_0)$, $t$ being in the direction of $\gamma$.
We will define
$s$ by the following: a point $(q,p)$ on the energy level of
 $\gamma$ will be assigned coordinates $(z,t,1)$ and the point 
$(q,sp)$ will be assigned coordinates $(z,t,s)$. It is clear that in 
an interval $s\in (a,b), a>0$, this  gives a system of coordinates.

It is interesting to notice that $(0,t,1)$ is a parametrization of 
$\gamma$, whereas the cylinder $(0,t,s)$ is foliated by circles
$s=c$ invariant
 under the flow $h_0^t$: each one corresponds to a reparametrization of 
$\gamma$, by rescaling the velocity by $s$.

The map $h_0^d$ leaves the cylinder invariant and in fact induces a 
monotone twist map on it :
$$
h_0^d(0,t,s)=(0,t + (s-1)d, s)
$$
Now, remember that the geodesic flow of a manifold with strictly
negative curvature
is Anosov. This translates into: in the subspace tangent at a point $(0,t,s)$
to the $z$ coordinate,
$Dh_0^d$ has no eigenvalue equal to 1
( we can assume the splitting tangent to $t,s$ to be invariant by $Dh_0^t$)
. Hence, in the $(z,t,s)$ coordinates:
$$
Dh_0^d(0,s,t)
\pmatrix{A& & \cr
& 1&d\cr
& 0&1\cr}
$$
where $A$ has no eigenvalue 1. Hence $Dh_0^d(0,t,s)$ has only the vector
tangent to $\gamma$ as eigenvector with eigenvalue 1 (I am endebted to
Leonid Polterovitch for giving me the idea of this argument).

We now have to prove that $G^0_{m,d}$ is a retract of $O_{m,d}$.
Define the following map
$$
\eqalign{\rho:O_{m,d}(2N) &\to O_{m,d}(N)\cr
             (q_0,q_1,q_2,\ldots, q_{2Nd}) &\to (q_0,q_2,\ldots,q_{2k}\ldots,
q_{2Nd})\cr }
$$
It is not hard to see that $\rho$ induces a diffeomorphism on $G^0_{m,d}$:
the projection on the 0th factor would itself give a diffeomorphism.
We claim that the image $\underline G$ of $G^0_{m,d}$ under $\rho$ is a deformation
retraction of the image $\underline O$ of $O_{m,d}$ in $O_{m,N}$. Call $r$ this
retraction. Then $\rho^{-1}_{|\underline G}\circ r \circ \rho $ is a retract of 
$O_{m,d}$ to $G^0_{m,d}$, as we want to prove.

We now construct the map $r$. Decompose $h_0^d=(h_0^{1\over N})^{Nd}$.
Since $h_0^{1\over N}$ is a 
symplectic twist map , we can rig
up the variational setting relative to this decomposition.
 Call 
$$\underline W(\qq)=\sum_{k=1}^{k=Nd} S(q_k,q_{k+1})
$$
where $S$ is the generating function of $h_0^{1\over N}$.
In this case, the isolating block $\underline B=\rho B$ is a repeller block
for the gradient flow of $\underline W$ (see Remark 5.7). Hence $\underline W$ has attains a 
minimum value, say $a$, in  the interior of $\underline  B$. It has to 
 be at a point 
in $\underline G$, which
contains all the  critical points of $\underline W$, as we remarked above for
$G^0_{m,d}$. Hence on all of $\underline G$
, $\underline W$ must equal $a$. Since
we can choose $\underline O$ to be exhausted by an increasing sequence of
repeller neighborhoods of the same type as $\underline P$, $a$ is actually a 
global minimum for the function $\underline W$ in $\underline O$.

 The same argument as for $G^0_{m,d}$ shows that 
$\underline G$ is normally hyperbolic. In particular, this implies that the set
$\underline W\leq a+\epsilon$ forms a tubular neighborhood of $\underline G$ 
([DNF],\S 20).
 Then a standard
argument ([Mi] Theorem 3.1)
in Morse theory shows that, since there are no other critical
points but those in the level $\underline W=a$, the set $\underline W\leq a + \epsilon $
 must be a deformation retraction of $\underline O$. Finally, $\underline G$ is
a deformation retraction of $\underline W \leq a + \epsilon$, since the latter is 
a tubular neighborhood of $\underline G$. This finishes the construction of
$r$. 

Finally, we indicate how all these features go through in the quotient by $\sigma$.

To check that $G^0_{m,d}/\sigma$ is normally hyperbolic, we just note
that this notion 
 is a local one, in the tangent space, and the quotient
map is a local diffeomorphism.
It can be checked that the above construction of the retraction map is
$\sigma$ invariant.
And, finally,
the quotient of  our invariant set $G^0_{m,d}$ 
(see 7.3) is
 a circle, as noted in 7.1.2.
 This concludes the proof of 
Lemma 7.2.2.\qed }

To finish the proof of Theorem 2, we use Floer's Lemma, as in the proof of 
Theorem 1, to find that there is an invariant set $G_{m,d}$ for the flow
$\z$ in $O_{m,d}/\sigma$ which is such that:
$$
H^*(G^0_{m,d}/\sigma)=H^*(\bbbs^1)\hookrightarrow H^*(G_{m,d})
$$
Since $cl(\bbbs^1)=sb(\bbbs^1)=2$, in all cases, we will get at least
 2 distinct {\it orbits} of type $(m,d)$.

\qed

\titlea{Appendix A}{ Linearized gradient flow vs.linearized Hamiltonian flow}

Suppose that $(q_0,p_0)=x_0$ is a fixed point for $F$. We want to solve the equation:
$$
DF_{x_0}(v) = \lambda  v \leqno (A.1)
$$
with $v \in T( T^*M)_{x_0}$. In terms of Hamiltonian flow
, we want to find the Floquet multipliers of the 
periodic orbit corresponding to
$x_0$. 

In the  $(q_k,q_{k+1})$ coordinates, we want to express a condition
on the orbit  $(\delta q_k, \delta q_{k+1})$ 
of a tangent vector
$(\delta q_1, \delta q_2)$
under the successive differentials of the maps $F_{k-1}$ 
 along the  given orbit. A way to do it is the 
following ([M-M]):
If $\overline q$ corresponds to the  orbit of $x_0$ under the  the successive 
$F_k$'s, it
must satisfy:
$$
{\partial W(\qq)\over \partial q_k}
=\partial_2 S_{k-1}(q_{k-1}, q_k) + \partial_1 S_k(q_k, q_{k+1})=0
$$
(see 4.3). Therefore, a ``tangent orbit'' $\delta \qq$ must satisfy:

$$
S_{21}^{k-1} \delta q_{k-1} + (S_{11}^k + S_{22}^{k-1})\delta q_k
+S_{12}^k \delta q_{k+1}=0 \leqno (A.2)
$$
where we have abbreviated: 
$$
S_{ij}^k=\partial_{ij}S_k(q_k, q_{k+1}).
$$

When $\qq$ 
corresponds to a fixed point $(q_0,p_0)$. 
Equation A.1 translates, in terms of the $\delta \qq$, to:
$$
\delta q_{2N}= \lambda \delta q_0 \leqno (A.3)
$$
Equations (A.2) and (A.3)
 can be put in matrix form as $M(\lambda)\delta \qq =0$ where $M(\lambda)$ 
is the following
$2Nn \times 2Nn$ tridiagonal matrix:
$$
M(\lambda)=\pmatrix{
S_{22}^0+S_{11}^1 & S_{12}^1 & 0 & \ldots & 0 & {1\over \lambda}S_{21}^0\cr
S_{21}^1 & S_{22}^1+S_{11}^2 & S_{12}^2 & \ddots &  & 0 \cr
0 & S_{12}^2 & & & &\vdots \cr
\vdots & \ddots &  & & & 0\cr
0 & \ldots & 0 & & & S_{12}^{2N-1}\cr
\lambda S_{12}^0 & 0 & \ldots & 0 & S_{21}^{2N-1}& S_{22}^{2N-1}+S_{11}^0\cr}
$$
Hence the eigenvalues of $DF_{x_0}$ are in one to one correspondence with
the values $\lambda$ for which $detM(\lambda)=0$. More precisely, to each
vector $v$ solution of (A.1) corresponds one and only one vector
$\delta \qq$ solution of $M(\lambda)\delta \qq =0$. Setting $\lambda =1$,
this proves Lemma 6.3.

\medskip {\bf Remark (A.4)}{ This construction  can be given a symplectic interpretation:
the Lagrangian manifolds  graph$(dW)$ and graph$(F)$ are related by symplectic
reduction. Lemma 6.3 can then be restated in terms of the invariance of
a certain Maslov index under reduction ([V]).

\titlea{}{Appendix B: Twist maps and linking of spheres}

In this appendix, we present the proof given in [BG] of
the original conjecture of Arnold in the restrictive case of 
symplectic twist maps (Theorem 3 in the introduction).

To that effect, we have to give our interpretation of what linking
of spheres in $\partial B^*M$ is.

Call $\Delta_q$ the fiber of $B^*M$ over $q$, and $\partial\DD$ its
boundary in $\partial B^*M$. Then $\partial \DD$ is an $n$ dimensional
sphere. It make sense to talk about its linking with its image
$F(\partial\DD)$ in $\partial B^*M$: the latter set has dimension
$2n-1$ and the dimensions of the spheres add up to $2n-2$.

 We first restrict ourselves to the 
case when the two spheres $\partial\DD$ and $F(\partial\DD)$ are in a trivializing
neighborhood in  $\partial B^*M$, say $U\times \bbbs^n=E$.

The type of linking of $F(\partial\DD)$ with $\partial\DD$ should
be given by the class $[F(\partial\DD)]\in H_{n-1}(\partial E\backslash
\partial \DD)$

More precisely, we have:
$$
\eqalign{H_{n-1}(\partial E\backslash \partial\DD) &\cong 
H_{n-1}\left(\bbbs^{n-1}\times (\bbbr^n-\{0\})\right)\cr
&\buildrel \rm Kunneth \over \cong
 H_{n-1}(\bbbr^n-\{0\})\oplus H_{n-1}(\bbbs^n)\cr}.\leqno (B.1)
$$
Thus, taking $\partial\DD$ from $\partial E$ creates a new generator in the 
$n-1$st homology, i.e. the generator $b$ of $H_{n-1}(\bbbr^n-\{0\})$.

\definition{(Linking condition)}{
We say that the spheres $F(\partial\DD)$ and  $\partial\DD$ link
in $\partial E$ if they do not intersect and if the decomposition
of $[F(\partial\DD)]$ in the direct sum  in (B.1) has a non zero term in
its $H_{n-1}(\bbbr^n-\{0\})$ factor. We will say that the symplectic twist map
$F$ satisfies the
linking condition if for all $q \in M$ these spheres link in  $\partial E$ for
some trivializing neighborhood $E$
\fonote{Here, as a convention, a trivializing neighborhood will always be 
homeomorphic to $\bbbb^n\times\bbbr^n$.}
}

If $F$ is a symplectic twist map, it turns out that this is a well defined
 characterization of linking: we can always  construct a trivializing
neighborhood containing both $\partial\DD$ and its image. Indeed, take
$T^*(\pi\circ F(\DD))$ (homeomorphic to $\bbbb^n\times \bbbr^n$
 since $F$ is twist)
if $q$ is in $\pi\circ F(\DD)$. If not join $q$
to this set by a path, and fatten this path. The union of the set and the
fattened path is homeomorphic to a ball. Hence the bundle over this ball
is trivial. 

Moreover, it turns out that if the spheres link in one trivializing 
neighborhood,
 they do in all of them, as a consequence
of the following

\lemma{B.2}{If F is a symplectic twist map, the following are equivalent:}
{\it
\item{a)} The spheres $\partial\DD$ and $F(\partial\DD)$ link in
some trivializing neighborhood in $\partial B^*M$
\item{b)} The fiber $\DD$ and its image $F(\DD)$ intersect
in one point of their interior.}

\medskip 
{\bf Remark B.3}{We can also define the linking condition for a map $F$ of $B^*M$
which is not necessarily symplectic twist.
If the covering space of $M$ is $\bbbr^n$, we say that $F$ satisfies the 
linking condition if at least one of its lifts does (the trivializing neighborhood
is taken to be $\overline M \times \bbbr^n \cong \bbbr^{2n}$ in this case.)
If $M$ is not covered by $\bbbr^n$, Lemma B.2 suggests that 
we may take as a linking condition that
the intersection number $\sharp (\DD \cap F(\DD))$ is $\pm 1$.

\proof{Suppose that $E$ is a trivializing neighborhood containing
the 2 spheres. 
We complete B.1 into  the following commutative diagram:
$$
\matrix{H_{n-1}(\partial E\backslash\partial \DD)& \cong &
 H_{n-1}(\bbbr^n-\{0\})\oplus H_{n-1}(\bbbs^n)\cr
\Big\downarrow\ i_* & & \Big\downarrow \ j_* \cr
H_{n-1}(E\backslash\DD)&\cong&H_{n-1}\left((\bbbr^n -\{0\}\times B^n)\right)
\cr}$$
where $i,j$ are inclusion maps. It is clear that $j_*b$ generates
$$
H_{n-1}\left((\bbbr^n-\{0\})\times B^n\right)\cong H_{n-1}\left((\bbbr^n-\{0\})
\times \bbbr^n\right).
$$
 The last group 
measures the 
(usual) linking number of a sphere with the fiber $\pi^{-1}(q)$
in $T^*E\cong \bbbr^{2n}$. But it is well known that such a number
is the intersection number of any ball bounded by the sphere with
the fiber $\pi^{-1}(q)$, counted with orientation. In our case, 
where the sphere considered is $F(\partial\DD)$
this number can only be 0 or 1 or -1, because of the twist condition (see
Remark 3.2).

Conversally, if $\DD$ and $F(\DD)$ intersect in their interior
, then their bounding
spheres must lie on the trivializing neighborhood over $F(\DD)$,
and must link. \qed }
\medskip 
{\bf Remark  B.4}{If all fibers intersect their image under a twist map, i.e.
if the linking condition is satisfied, then the intersection number
must be uniformally 1 or -1: we could call $F$ a positive twist
map in the first case, a negative twist map in the second case. Of
course, this corresponds to the same notion in dimension 2.}
\medskip
We can now prove Theorem 3.

>From Lemma 4.3, fixed points of F correspond to critical points
of $q\to S(q,q)$. This function only make sense for all $q$ in $M$ if the
diagonal in $M\times M$ is in the image of $U$ by the embedding
$\psi$ (see Definition 3.1). This is exactly the case when $q \in F(\DD)$
for all $q$, i.e., from Lemma B.2, exactly when the linking condition is
satisfied. Hence $F$ has as many fixed points as the function
$q\to S(q,q)$ has critical points on $M$.\qed

\medskip
To our knowledge, Arnold's original conjecture is still open, even in the case
$M=\bbbt^n$.

\begrefchapter{References}

\refno{[Ar1]}V. Arnold, `` Sur une propri\'et\'e topologique des applications
globalement cannoniques de la m\'ecanique classique'', C.R. Acad.Sc. Paris,
t.261 (1965). Groupe 1.

\refno{[Ar1]} V.I. Arnold: ``Mathematical Methods of Classical Mechanics''
 (Appendix 
 9), Springer-Verlag 1978.

\refno{[AL]} S. Aubry and P.Y. LeDaeron: ``The discrete Frenkel-Kontarova
model and its extensions I. Exact results for ground states'', {\it Physica} 8D
 (1983), 381-422.

\refno{[BG]}A.Banyaga and C. Gol\'e : ``A
remark on  a conjecture of Arnold: linked 
spheres and fixed points'', to appear in proc. conf. on Hamiltonian systems
and celestial mechanics, Guanajuato

\refno{[BK]} D. Bernstein and A.B. Katok: ``Birkhoff periodic orbits for small
perturbations of completely integrable Hamiltonian systems with convex
Hamiltonians'', {\it Invent. Math.} {\bf 88} (1987), 225-241.

\refno{[BP]}M. Bialy and L. Polterovitch,
`` Hamiltonian diffeomorphisms and Lagrangian distributions'', preprint,
Tel Aviv University, (1991)

\refno{[Ch1]} M. Chaperon: ``Quelques questions de G\'eometrie Symplectique'',
S\'eminaire Bourbaki no. 610 (1982/83).

\refno{[Ch2]} M. Chaperon, ``Une id\'ee du type ``g\'eod\'esiques bris\'ees''
pour les syst\`emes hamiltoniens'', {\it C.R. Acad. Sc., Paris}, 298, S\'erie I,
no 13, 1984, p 293-296.

\refno{[C]} Conley, C.C.: ``Isolated invariant sets and the Morse index''
, CBMS,
Regional Conf. Series in Math., Vol.38 (1978).

\refno{[CZ1]} C.C. Conley and E. Zehnder: ``The Birkhoff-Lewis 
fixed point theorem
 and a conjecture of V.I. Arnold'', {\it Invent. Math.} (1983).

\refno{[CZ2]} C.C. Conley and E. Zehnder: ``Morse type index theory for 
Hamiltonian equations'', {\it Comm. Pure and Appl. Math.}
 XXXVII (1984), pp. 207--253.

\refno{[D]} R. Douady: ``Stabilit\'e ou instabilit\'e des points fixes 
elliptiques'', Ann. Sci. Ec. Norm. Sup. 4\`eme s\'erie, t.21 (1988), pp. 1--46

\refno{[DNF]} B. Doubrovine, S. Novikov, A. Fomenko : ``G\'eom\'etrie
contemporaine'', vol 3, Editions Mir, Moscow, (1987) (see also english
translation in Springer--Verlag)

\refno{[F1]} A. Floer: ``A refinement of Conley index and an application to the
stability
of hyperbolic invariant sets'',
 Ergod. Th. and Dyn. Sys., vol 7, ( 1987).

\refno{[F2]} A. Floer: ``Morse theory for Lagrangian intersections'',
J. Diff. Geom. 28 (1988) pp. 513--547

\refno{[G1]} C. Gol\'e, ``Periodic points for monotone symplectomorphisms of
$\Tn \times \Rn$''; Ph.D. Thesis, Boston University (1989)

\refno{[G2]} C. Gol\'e, ``Ghost circles for twist maps''; IMA preprint (1990), to appear in
J. of Diff. eq.

\refno {[GH]} C. Gol\'e and G.R. Hall:
``Poincar\'e's proof of Poincar\'e's last geometric theorem''
To appear in the proceedings for the workshop on twist maps, IMA
,Springer

\refno{[GHL]}S. Gallot, D. Hulin and J. Lafontaine: ``Riemannian Geometry'',
Springer-Verlag, (1987).

\refno{[Gr]} M.J. Greenberg: ``Lectures on Algebraic Topology''
, Math. lecture note
 series (5th printing, 1977).

\refno{[H]} M.R. Herman: ``Existence et non existence de Tores invariants
par des diff\'eomorphismes symplectiques'', S\'eminaire sur les 
\'equations aux d\'eriv\'es partielles, Ecole Polytechique, Palaiseau (1988).

\refno{[J]} F.W. Josellis: ``Global periodic orbits for Hamiltonian systems
on ${\bf T}^n\times {\bf R}^n$''
Ph.D. Thesis Nr. 9518, ETH Z\"urich, (1991).

\refno{[Ka]} A. Katok: ``Some remarks on Birkhoff and Mather twist map 
theorems'',
 Ergod. Th. \& Dynam. Sys. (1982), 2, 185-194.

\refno{[Kl]}W. Klingenberg: ``Riemannian geometry'', de Gruyter studies
in Mathematics (1982).

\refno{[K-M]} H. Kook and J. Meiss: ``Periodic orbits for Reversible,
 Symplectic
 Mappings'', Physica D 35 (1989) 65--86.

\refno{[L]} P. LeCalvez: ``Existence d'orbites de Birkhoff G\'en\'eralis\'ees
pour les diff\'eo  morphismes conservatifs de l'anneau'', Preprint, Universit\'e
Paris-Sud, Orsay.

\refno{[McD]} D. McDuff:  ``Elliptic methods in symplectic geometry'',
Bull. Am. Math. Soc., vol 23, n.2, (1990)

\refno{[MM]]} R.S. Mackay and J. Meiss: ``Linear stability of periodic orbits
 in Lagrangian systems'', {\it Phys. Lett.} {\bf 98A}, 92 (1983).

\refno{[MMS]}R.S McKay, J. Meiss and J.Stark: ``Converse KAM theory
for symplectic twist maps'', Nonlinearity 2 (1989) pp 469--512.

\refno{[Ma]} J. Mather: `` Action minimizing invariant measures for
positive definite Lagrangian systems'', preprint, ETH Z\"urich, (1989)

\refno{[Mi,2]} J. Milnor: ``Morse Theory'', Princeton University Press.

\refno{[Mo1]} J.Moser: ``Proof of a generalized form of a fixed point
theorem due to G.D. Birkhoff'', Lecture Notes in Mathematics, Vol. 597:
Geometry and Topology, pp. 464-494. Springer (1977)

\refno{[Mo2]} J. Moser: ``Monotone twist mappings and the calculus of
 variations'',
Ergod. Th. and Dyn. Sys., Vol 6 (1986), p 401-413

\refno{[V]} C. Viterbo: `` Intersection de sous vari\'et\'es Lagrangiennes, 
fonctionnelles d'action et indice de syst\`emes Hamiltoniens'', Bull. Soc. math.
France 115, (1987), pp. 361--390
\endref	

\bye